\documentclass[a4paper,12pt,leqno]{amsart}

\usepackage{amssymb,latexsym}
\usepackage{mltex}
\usepackage{amsmath}
\usepackage{epsfig}
\usepackage{graphics}
\usepackage{pst-3d}
\usepackage{pb-diagram}
 \usepackage{lscape} 
\usepackage{enumerate}

\newtheorem{defn}{Definition}[section]
\newtheorem{thm}{Theorem}   
     
\newtheorem{lem}[defn]{Lemma}     
\newtheorem{prop}[defn]{Proposition}

\newcommand{\Mndr}{\mathcal{M}(n,\delta,R)}
\newcommand{\nablab}{\overline{\nabla}}
\newcommand{\expo}{{\exp_{p_0}}}
\newcommand{\Ric}{\mathrm{Ric}}
\newcommand{\diff}{\mathrm{d}\,}
\newcommand{\supp}{\mathrm{supp}}
\newcommand{\td}{t_{\delta}}
\newcommand{\sd}{s_{\delta}}
\newcommand{\cd}{c_{\delta}}
\newcommand{\lgra}{\longrightarrow}
\newcommand{\Hinf}{||H||_{\infty}}

\newcommand{\HHinf}{||H||^2_{\infty}}
\newcommand{\Binf}{||B||_{\infty}}
\newcommand{\Hp}{||H||_{2p}}
\newcommand{\HHp}{||H||^2_{2p}}

\newcommand{\diam}{\mathrm{diam}}

\newcommand{\ddiv}{\mathrm{div}\,}

\newcommand{\iinf}{\mathrm{inf}\,}

\newcommand{\trace}{\mathrm{tr\,}}

\newcommand{\lto}{\ensuremath{\longrightarrow}}

\newcommand{\Ss}{\mathbb{S}}
\newcommand{\HH}{\mathbb{H}}

\newcommand{\R}{\mathbb{R}}

\newcommand{\Md}{\mathbb{M}^{n+1}(\delta)}

\newcommand{\function}[5]
{\begin{eqnarray*}\begin{array}{r@{}ccl}
 #1\;\colon\;  & #2 &\lto & #3 \\[.05cm]  
  & #4 &\longmapsto  & #5 
\end{array}\end{eqnarray*}
}
\newcommand{\beqt}{\begin{equation}}  \newcommand{\eeqt}{\end{equation}}
\newcommand{\bal}{\begin{align}}      \newcommand{\eal}{\end{align}}
\newcommand{\ba}{\begin{array}}      \newcommand{\ea}{\end{array}}
\newcommand{\bc}{\begin{center}}     \newcommand{\ec}{\end{center}}
\newcommand{\be}{\begin{enumerate}}  \newcommand{\ee}{\end{enumerate}}
\newcommand{\beq}{\begin{eqnarray}}  \newcommand{\eeq}{\end{eqnarray}}
\newcommand{\beQ}{\begin{eqnarray*}} \newcommand{\eeQ}{\end{eqnarray*}}
\newcommand{\bi}{\begin{itemize}}    \newcommand{\ei}{\end{itemize}}
\newcommand{\bt}{\begin{tabular}}    \newcommand{\et}{\end{tabular}}
\newcommand{\finpreuve}{\hfill\square\\}

\title[Extrinsic radius pinching]{Extrinsic radius pinching for hypersurfaces of space forms}
\author{Julien ROTH}
\date{\today}

\begin{document}
\maketitle
\begin{center}
Institut \'Elie Cartan, Universit\'e Henri Poincar\'e, Nancy I, B.P. 239
\end{center}
\begin{center}
54506 Vand\oe uvre-L\`es-Nancy Cedex, France
\end{center}
\begin{center}
roth@iecn.u-nancy.fr
\end{center}
\begin{abstract}
We prove some pinching results for the extrinsic radius of compact hypersurfaces in space forms. In the hyperbolic space, we show that if the volume of $M$ is $1$, then there exists a constant $C$ depending on the dimension of $M$ and the $L^{\infty}$-norm of the second fundamental form $B$ such that the pinching condition $\tanh(R)<\frac{1}{\Hinf}+C$ (where $H$ is the mean curvature) implies that $M$ is diffeomorphic to an $n$-dimensional sphere. We prove the corresponding result for hypersurfaces of the Euclidean space and the sphere with the $L^p$-norm of $H$, $p\geq2$, instead of the $L^{\infty}$-norm.\\\\
\end{abstract}

\section{Introduction}
\label{introduction}
Let $(M^n,g)$ be a compact, connected and oriented $n$-dimensional Riemannian manifold without boundary isometrically immersed by $\phi$ into the $(n+1)$-dimensional simply connected space-form $(\Md,g_{can})$ of sectional curvature $\delta$ with $n\geq2$. First, let us recall the definition of the extrinsic radius of M.
\begin{defn}
The extrinsic radius of $(M,g)$ is the number
$$R=R(M)=\iinf\left\lbrace r>0\big|\ \exists x\in\Md\ \text{s.t.}\ \phi(M)\subset B(x,r)\right\rbrace,$$
where $B(x,r)$ is the open ball of center $x$ and radius $r$ in $\Md$. 
\end{defn}
Throughout this paper, we denote respectively by $B(x,r)$, $\overline{B}(x,r)$ and $S(x,r)$ the open ball, the closed ball and the sphere of center $x$ and radius $r$ in $\Md$. An immediate consequence of the above definition is that there exists $p_0\in\Md$ such that $\phi(M)\subset \overline{B}(p_0,R)$ and $\phi(M)\cap S(p_0,R)\neq\emptyset$. Moreover, it is a well-known fact that the extrinsic radius is bounded from below in terms of the mean curvature. More precisely, we have the following estimate due to Hasanis and Koutroufiotis (\cite{HK}) for $\delta=0$ and Baikoussis and Koufogiorgos (\cite{BK}) for any $\delta$
\beqt\label{lowerboundh}
\td(R)\geq\frac{1}{\Hinf},
\eeqt
where $\td$ is the function defined in Section \ref{preliminaries} and $H$ the mean curvature of the immersion. Note that for $\delta>0$, the image $\phi(M)$ is assumed to be contained in a ball of radius less than $\frac{\pi}{2\sqrt{\delta}}$, that is an open hemisphere.
Moreover, equality in (\ref{lowerboundh}) is characterized by geodesic hyperspheres.\\
\indent
A natural question is the following: Is there a constant $C$, depending on a minimal number of geometric invariants, such that if we have the pinching condition
\beqt
\td(R)<\frac{1}{\Hinf}+C,\tag{$P_C$}
\eeqt
then $M$ is closed, in a certain sense, to a sphere?\\
\indent
Many pinching results are known for geometric invariants defined on Riemannian manifolds with positive Ricci curvature, as the intrinsic diameter (\cite{Esch,Ili,Wu}), the volume, the radius (\cite{Col1,Col2}) or the intrinsic lower bound of Lichnerowicz-Obata of the first nonzero eigenvalue of the Laplacian in terms of lower bounds of the Ricci curvature (\cite{Cro,Ili,Pet}).\\
\indent
For instance, concerning the intrinsic diameter, under the hypothesis that $(M^n,g)$ is a complete Riemannian manifold with $\Ric\geq n-1$, Myers gave the well-known upper bound
$$\diam(M^n,g)\leq\diam(\Ss^n, can)=\pi.$$
In particular, $M$ is a compact manifold.\\
\indent
S. Ilias proved in \cite{Ili} that there exists an $\varepsilon$ depending on $n$ and an upper bound of the sectional curvature so that if $\Ric\geq n-1$ and $ \diam(M)>\pi-\varepsilon$, then $M$ is homeomorphic to $\Ss^n$.\\
\indent
Petersen and Sprouse gave in \cite{PS} a generalization of the Theorem of Myers with a less restrictive assumption on the Ricci curvature. They assume that $\Ric$ is almost bounded from below by $n-1$ in an $L^p$-sense. Then under this hypothesis, $\diam(M^n,g)\leq\pi+\varepsilon$.\\
\indent
With a similar hypothesis on the Ricci curvature, E. Aubry (\cite{Aub}, Theorem 5.24) proved that if $\diam(M^n,g)\geq\pi-\varepsilon$ for $\varepsilon$ small enough depending on an upper bound of the sectional curvature, then $M^n$ is homeomorphic to $\Ss^n$.\\
\indent
In this paper, the hypothesis on the Ricci curvature is replaced by the fact that $M$ is isometrically immersed in a standard space form. Moreover, as we will see, the upper bound of the sectional curvature will be replaced by the $L^{\infty}$-norm of the mean curvature or that of the second fundamental form. Recently, under the hypothesis that $M$ is isometrically immersed in the Euclidean space, Colbois and Grosjean (see \cite{CG}) proved a pinching result on the first eigenvalue of the Laplacian. More precisely, they proved that there exists a constant $C$ depending on $n$ and the $L^{\infty}$-norm of the second fundamental form such that if $\frac{n}{V(M)^{1/p}}||H||_{2p}^2-C<\lambda_1(M)$, then $M$ is diffeomorphic to an $n$-dimensional sphere.\\
\indent
We keep on with studying hypersurfaces where little is known about pinching results. Indeed, we give pinching results for the extrinsic radius, which is the extrinsic analogue to the diameter, for hypersurfaces of the Euclidean space and hypersurfaces of the sphere and the hyperbolic space too.\\
\indent
For more convenience, we denote by $\Mndr$ {\it the set of all compact, connected and oriented $n$-dimensional Riemannian manifolds without boundary  isometrically immersed into $\Md$ of extrinsic radius $R$ and volume equal to $1$}. In the case $\delta>0$, we assume that $M$ lies in an open hemisphere of $\Ss^{n+1}(\delta)$.
\begin{thm}\label{thm1}
Let $(M^n,g)\in\Mndr$ and let $p_0$ be the center of the ball of radius $R$ containing $M$. Then for any $\varepsilon>0$, there exists a constant $C_{\varepsilon}$ depending only on $n$, $\delta$ and the $L^{\infty}$-norm of the mean curvature such that if 
\beqt
\td(R)<\frac{1}{\Hinf}+C_{\varepsilon}\tag{$P_{C_{\varepsilon}}$}
\eeqt
then 
\begin{enumerate}[i)]
\item $\phi(M)\subset \overline{B}\big(p_0,R\big)\setminus B\big(p_0,R-\varepsilon\big).$
\item $\forall x\in S\big(p_0,R\big),\quad B(x,\varepsilon)\cap\phi(M)\neq\emptyset$.
\end{enumerate}
\end{thm}

\noindent
{\bfseries Remark.} We will see in the proof that $ C_{\varepsilon}\lgra0$ when $\Hinf\lgra\infty$ or $\varepsilon\lgra0$.\\\\
\indent
We recall that the Haussdorff-distance between two compact subsets $A$ and $B$ of a metric space is given by
$$d_H(A,B)=\iinf\left\lbrace \eta\big|\;B\subset V_{\eta}(A)\quad \text{and}\quad A\subset V_{\eta}(B)\right\rbrace $$
where for any subset $A$, $V_{\eta}(A)$ is the tubular neighborhood of $A$ defined by $V_{\eta}(A)=\left\lbrace x\big|\ {\rm dist}(x,A)<\eta\right\rbrace $. So the points $i)$ and $ii)$ of Theorems \ref{thm1} imply that 
$$d_H\left(M,S(p_0,R)\right)\leq\varepsilon.$$
\indent 
If the pinching condition is strong enough, with a control on the $L^{\infty}$-norm of the second fundamental form instead of the $L^{\infty}$-norm of the mean curvature, we obtain that $M$ is diffeomorphic to a sphere and almost isometric to a geodesic sphere in the following sense:
\begin{thm}\label{thm2}
Let $(M^n,g)\in\Mndr$ and let $p_0$ be the center of the ball of radius $R$ containing $M$. Then there exists a constant $C$ depending only on $n$, $\delta$ and the $L^{\infty}$-norm of the second fundamental form such that if $(P_C)$ is true, then $M$ is diffeomorphic to $S(p_0,R)$.\\
\indent
More preciesly, there exists a diffeomorphism $F$ from $M$ into the geodesic hypersphere $S(p_0,R)$ of radius $R$ which is a quasi-isometry. That is, for all $\theta\in]0,1[$, there exists a constant $C$ depending on $n$, $\delta$, $\Binf$ and $\theta$ such that the pinching condition $(P_C)$ implies
$$\big||dF_x(u)|^2-1\big|\leq\theta,$$
for all unit vector $u\in T_xM$.
\end{thm}
\noindent
{\bfseries Remark.}
In the two above Theorems, we assume that $V(M,g)=1$. By homothety, we can deduce the same results for manifolds with arbitrary volume. Indeed, $(M,g')\in\mathcal{M}(n,\delta',R')$, with $g'=V(M)^{-2/n}g$, $\delta'= V(M)^{2/n}\delta$ and $R'=V(M)^{-1/n}R$.
\\\\
\indent
We will see in Section \ref{preliminaries} that in the case $\delta\geq0$, Inequality (\ref{lowerboundh}) can be improved by replacing $\Hinf$ by $\Hp$ (see Proposition \ref{lowerboundhp}). Moreover, the equality is also caracterized by geodesic hyperspheres. Therefore, we can consider the corresponding pinching problem. Theorems \ref{thm3} and \ref{thm4} give the analogue of Theorems \ref{thm1} and \ref{thm2} for this integral lower bound.\\\\
\noindent
{\bfseries Acknowledgement.} The author would like to thank Jean-Fran\c{c}ois Grosjean and Oussama Hijazi for their support.

\section{Preliminaries}
\label{preliminaries}
First, let us introduce the following functions:
$$\sd(t)=\left\lbrace \begin{array}{ll}
\frac{1}{\sqrt{\delta}}\sin( \sqrt{\delta}\,t) & \text{if}\ \delta>0\\
t & \text{if}\ \delta=0 \\ 
 \frac{1}{\sqrt{-\delta}}\sinh( \sqrt{-\delta}\,t) & \text{if}\ \delta<0
\end{array}\right. $$
and 
$$\cd(t)=\left\lbrace \begin{array}{ll}
\cos(\sqrt{\delta}\,t) & \text{if}\ \delta>0\\
1 & \text{if}\ \delta=0 \\ 
\cosh( \sqrt{-\delta}\,t) & \text{if}\ \delta<0
\end{array}\right. $$
\indent 
We can easily check that $\cd^2+\delta\sd^2=1$, $\sd'=\cd$ and $\cd'=-\delta\sd$. Moreover, we define the function $\td=\dfrac{\sd}{\cd}$ which satisfies $\td'=1+\delta\td^2$.\\
\indent
Throughout this paper, we consider a  Riemannian manifold $(M^n,g)$ of $\Mndr$. The second fundamental form $B$ of the immersion is defined by
$$B(X,Y)=\left\langle \nablab_X\nu,Y\right\rangle ,$$
where $<\cdot,\cdot>$ and $\nablab$ are respectively the Riemannian metric and the Riemannian connection of $\Md$. The mean curvature of the immersion is $$H=\frac{1}{n}\trace(B).$$\\
\indent
For any $p_0\in\Md$ let $\expo$ be the exponential map at this point. We consider $(x_i)_{1\leq i\leq n+1}$ the normal coordinates of $\Md$ centered at $p_0$. For $x\in\Md$, we denote by $r(x)=d(p_0,x)$ the geodesic distance from $p_0$ to $x$ on $(\Md,g_{can})$.\\
\indent
In what follows, $\nabla$ and $\nablab$ will be respectively the gradients associated to $(M,g)$ and $(\Md,g_{can})$. The corresponding Laplacians are $\Delta$ and $\overline{\Delta}$. The coordinates of $Z:=\sd(r)\nablab r$ in the normal frame are $\left( \frac{\sd(r)}{r}x_i\right)_{1\leq i\leq n+1}$. We denote by $X^T$ the projection of a vector field $X$  on the tangent bundle of $\phi(M)$.\\
\indent
Now let's recall some properties of the exponential map. First, $\expo$ is a radial isometry, {\it i.e.}, for each $x\in\Md$, we have
\beqt\label{radialisometry}
\left\langle ( \diff\expo)_X(X),(\diff\expo)_X(v)\right\rangle_x=\left\langle X,v\right\rangle_{p_0},
\eeqt
where $X=\expo^{-1}(x)$ and $v\in T_{p_0}\Md$. On the other hand, we have the following equalities (see Corollary 2.8 and Lemma 2.9 p153 in \cite{Sak}). If $v$ is a vector of $T_x\Md$ orthogonal to $\nablab r$, we have
\beqt\label{expequ1}
\big|\left(\diff\expo^{-1} \right)_{|x}(v) \big|^2_{p_0}=\frac{r^2}{\sd^2(r)}|v|^2_x,
\eeqt
and 
\beqt\label{expequ2}
\left\langle \nablab_v\nablab r,v\right\rangle=\frac{\cd(r)}{\sd(r)}|v|^2.
\eeqt
Moreover, $\nablab r$ is in the kernel of $\nablab\diff r$. In particular, for any \\$v\in T_x\Md$,
\beqt\label{nullspace}
\left\langle \nablab_{\nablab r}\nablab r,v\right\rangle=\left\langle \nablab_v\nablab r,\nablab r\right\rangle =0.
\eeqt
Finally, we give the following lemma (see \cite{Hei} or \cite{Gr2} for a proof):
\begin{lem}\label{lem1}
\begin{enumerate}[i)]
\item $\ddiv(Z^T)=n\cd(r)+nH\left\langle Z,\nu\right\rangle,$\\
\item $\delta\displaystyle\int_Mg(Z^T,Z^T)dv_g\geq n\int_M\Big( \cd^2(r)-|H|\cd(r)\sd(r)\Big) dv_g,$
\end{enumerate}
\end{lem}
\noindent
{\bfseries Remark.}
Note that, after integration, the first point in the case $\delta=0$ is nothing else but the Hsiung-Minkowski formula (see \cite{Hsi}).
\\\\
\indent
From this Lemma, we deduce the following estimates for the extrinsic radius in the case $\delta\geq0$.
\begin{prop}\label{lowerboundhp}
Let $(M^n,g)$ be a compact, connected and oriented $n$-dimensional Riemannian manifold without boundary isometrically immersed by $\phi$ into $\R^{n+1}$ or an open hemisphere of $\Ss^{n+1}(\delta)$. Then the extrinsic radius $R$ of $M$ satisfies
$$\td(R)\geq\frac{V(M)^{1/p}}{||H||_p},$$
for any $p\geq1$.
Moreover, equality holds if and only if $(M^n,g)$ is a geodesic hypersphere of radius $R$.
\end{prop}
{\it Proof:}
After integration, the first point of Lemma \ref{lem1} gives
$$\int_M\cd(r)dv_g\leq\int_M|H|\sd(r)dv_g.$$
Since $\sd$ is an increasing function and $\cd$ is a decreasing function, we get
$$\cd(R)V(M)\leq\sd(R)||H||_1.$$
The H\"older inequality gives the result for any $L^p$-norm.
Obviously, if $(M^n,g)$ is a geodesic hypersphere of radius $R$, then we have equality. Conversly, if equality holds, then 
$$\int_M\cd(r)dv_g=\cd(R)V(M),$$
which implies that $r\equiv R$ and $(M^n,g)$ is a geodesic hypersphere of radius $R$.
$\finpreuve$\\
We have the following pinching results corresponding to this inequality.
\begin{thm}\label{thm3}
Let $(M^n,g)\in\Mndr$ with $\delta\geq0$ and let $p_0$ be the center of the ball of radius $R$ containing $M$. Let $p\geq1$. Then for any $\varepsilon>0$, there exists a constant $C_{\varepsilon}$ depending only on $n$, $\delta$ and the $L^{\infty}$-norm of the mean curvature such that if 
\beqt
\td(R)<\frac{1}{\Hp}+C_{\varepsilon}\tag{$\widetilde{P}_{C_{\varepsilon}}$}
\eeqt
then 
\begin{enumerate}[i)]
\item $\phi(M)\subset \overline{B}\big(p_0,R\big)\setminus B\big(p_0,R-\varepsilon\big).$
\item $\forall x\in S\big(p_0,R\big),\quad B(x,\varepsilon)\cap\phi(M)\neq0.$
\end{enumerate}
\end{thm}

\begin{thm}\label{thm4}
Let $(M^n,g)\in\Mndr$ with $\delta\geq0$ and let $p_0$ be the center of the ball of radius $R$ containing $M$. Let $p\geq1$. Then there exists a constant $C$ depending only on $n$, $\delta$ and the $L^{\infty}$-norm of the second fundamental form such that if
\beqt
\td(R)<\frac{1}{\Hp}+C\tag{$\widetilde{P}_C$}
\eeqt
then $M$ is diffeomorphic and quasi-isometric to $S(p_0,R)$ in the sense of Theorem \ref{thm2}.
\end{thm}
\noindent
{\bfseries Remark.} In the case $\delta\geq0$, Theorems \ref{thm1} and \ref{thm2} are just corollaries of the two above theorems since if $V(M)=1$, $\Hinf\geq\Hp$.

\section{An $L^2$-approach to pinching}
\label{approach}
A first step in the proof of the pinching results is to prove that the pinching condition $(P_C)$
in the three cases, or $(\widetilde{P}_C)$, in the Euclidean or spherical case
implies that $M$ is close to a hypersphere in an $L^2$-sense.\\
\indent
For this, let's consider the functions $\varphi$ and $\psi$ definied by
$$\varphi(r)=\left\lbrace \begin{array}{l}
\td^2(R)-\td^2(r)\quad\text{if}\quad \delta<0, \\ \\
\sd^2(R)-\sd^2(r)\quad\text{if}\quad \delta\geq0
\end{array}\right. $$
and $$\psi(r)=\cd(r)|Z^T|.$$
\subsection{The Hyperbolic case}
\label{hyperboliccase}
In this section, we suppose $\delta<0$. Note that if the pinching constant $C$ satisfies
\beqt\label{pinchconst}
 C\leq\frac{1}{2}\left(\frac{1}{\sqrt{-\delta}}-\frac{1}{\Hinf} \right) =\alpha(\Hinf)
 \eeqt
then $(P_C)$ implies $\td(R)\leq \frac{1}{2}\left(\frac{1}{\sqrt{-\delta}}+\frac{1}{\Hinf} \right)=\beta(\Hinf)<\frac{1}{\sqrt{-\delta}}$ and $R$ is bounded from above by a constant depending only on $\Hinf$. In what follows, we assume that the pinching constant $C$ satisfies the relation (\ref{pinchconst}). We prove the following lemma
\begin{lem}\label{lemphi}
The pinching condition $(P_C)$ with $C\leq\alpha(\Hinf)$ implies
$$||\varphi||^2_2\leq A_1C,$$
where $A_1$ is a positive constant depending only on $\delta$ and $\Hinf$.
\end{lem}
{\it Proof:} Since $\td$ is an increasing function, we have
$$||\varphi||_2^2\leq\td^2(R)\int_M\left( \td^2(R)-\td^2(r)\right).$$  
Since $\td^2(R)-\td^2(r)\geq0$ and $\cd(r)\geq1$, so
\beQ
\int_M\left( \td^2(R)-\td^2(r)\right) &\leq&\td^2(R)\int_M\cd^2(r)-\int_M\sd^2(r)\\
&\leq&\td^2(R)\int_M\cd^2(r)-\frac{1}{\HHinf}\int_MH^2\sd^2(r).
\eeQ
Using the H\"older inequality, we get
\beQ
\int_M\left( \td^2(R)-\td^2(r)\right) &\leq&\td^2(R)\int_M\cd^2(r)-\frac{1}{\HHinf}\frac{\left( \int_M|H|\sd(r)\cd(r)\right)^2 }{ \int_M\cd^2(r)}
\eeQ
Now using the relation $ii)$ of Lemma \ref{lem1} and applying the pinching condition $(P_C)$ with $C$ satisfying (\ref{pinchconst}), we find
\beQ
\int_M\left( \td^2(R)-\td^2(r)\right)&\leq&\left( \td^2(R)-\frac{1}{\HHinf}\right)\int_M\cd^2(r) \\
 &\leq&\left( C^2+\frac{2C}{\Hinf}\right)\cd^2(R)
\eeQ
and
\beQ
||\varphi||^2_2&\leq&\sd^2(R)\left( C^2+\frac{2C}{\Hinf}\right)\leq A_1C,
\eeQ
where $A_1$ depends only on $n$, $\delta$ and $\Hinf$.$\finpreuve$\\
\indent
The next Lemma gives an upper bound for $||\psi||_2$ under the pinching condition.
\begin{lem}\label{lempsi}
The pinching condition $(P_C)$ with $C\leq\alpha(\Hinf)$ implies
$$||\psi||_2^2\leq A_2C+A_3||\varphi||_{\infty},$$
where $A_2$ depends only on $\delta$ and $A_3$ depends on $\delta$ and $\Hinf$.
\end{lem}
{\it Proof:}  First,  we  observe  that  $|Z^T|^2=|Z|^2-\left\langle Z,\nu\right\rangle^2$.  Since  \\$|Z|=\sd(r)$, we have
\beQ
||\psi||_2^2&\leq&\cd^2(R)\left[ \sd^2(R)-\frac{1}{\HHinf}\int_M\left( H^2\left\langle Z,\nu\right\rangle^2 \right) \right] .
\eeQ
Using H\"older inequality and Lemma \ref{lem1} $i)$, we get
\beQ
||\psi||_2^2&\leq&\cd^2(R)\left[ \sd^2(R)-\frac{1}{\HHinf}\left( \int_M\cd(r)\right)^2\right].\\
&=&\cd^2(R) \Bigg[\sd^2(R)-\frac{1}{\HHinf}\cd^2(R)-\frac{2\cd(R)}{\HHinf}\int_M\left(\cd(r)-\cd(R) \right) \\
&&-\frac{1}{\HHinf}\left( \int_M\cd(r)-\cd(R)\right)^2\Bigg] \\
&\leq&\left( \frac{2C}{\Hinf}+C^2\right)\cd^4(R)+K\big|\cd(r)-\cd(R)\big| 
\eeQ
where $K$ depends on $\delta$ and $\Hinf$. Since $\td(R)\leq\beta(\Hinf)<\frac{1}{\sqrt{-\delta}}$, we deduce that there exists $K'>0$ depending on $n$, $\delta$ and $\Hinf$ so that $\big|\cd(R)-\cd(r)\big|\leq K'||\varphi||_{\infty}$. This completes the proof.
$\finpreuve$\\
\subsection{The Euclidean and Spherical cases}
\label{sphericalcase}
Here, $\delta\geq0$ and if $\delta>0$ then we assume that $\phi(M)$ is contained in an open hemisphere ({\it i.e.} an open ball of radius $\frac{\pi}{2\sqrt{\delta}}$). First, note that the pinching condition $(\widetilde{P}_C)$ with $C<1$ and the fact that $V(M)=1$ imply that there exist two constants $\alpha_n$ and $\beta_n$ depending only on $n$ so that $\Hinf\geq\Hp\geq\alpha_n$ and $\td(R)\leq\beta_n$. Consequently, $R$ is bounded from above by a constant $\gamma_n$. That is an immediate consequence of the Sobolev following inequality due to Hoffman and Spruck (cf \cite{HS1} and \cite{HS2}) for a nonnegative fonction $f$, by taking $f\equiv1$
$$\left(\int_Mf^{n/(n-1)}dv_g\right)^{2(n-1)/n}\leq K_n\left(\int_M|H|f+|\nabla f|dv_g\right).$$
\indent
Note that this inequality is true without further assumptions for $\delta=0$. For $\delta>0$, some conditions on the sectional curvature and the injectivity radius $i\big(\Ss^{n+1}(\delta)\big)$ of $\Ss^{n+1}(\delta)$ and on the support of the function $f$ are needed. The first condition, $i\big(\Ss^{n+1}(\delta)\big)\geq\pi\delta^{-1}$, is satisfied since for the sphere $\Ss^{n+1}(\delta)$, we have $i\big(\Ss^{n+1}(\delta)\big)=\pi\delta^{-1}$. The second condition is $V(\supp(f))\leq(1-\alpha)\omega_n\delta^{-n}$, for some $0<\alpha<1$ and where $\omega_n$ is the volume of the $n$-dimensional Euclidean ball. This condition is automatically satisfied if $\phi(M)$ lies in an open hemisphere.\\
\indent
In the sequel, we assume that the pinching constant satisfies $C<1$. We have the following lemma
\begin{lem}\label{lemphi2}
The pinching condition $(\widetilde{P}_C)$ with $C<1$ implies
$$||\varphi||^2_2\leq \widetilde{A_1}C,$$
where $\widetilde{A_1}$ is a positive constant depending only on $n$ and $\delta$.
\end{lem}
{\it Proof:} Since $\sd$ is an increasing function and $\cd$ is a decreasing function, we have
\beQ
||\varphi||_2^2&\leq&\sd^2(R)\int_M\left( \sd^2(R)-\sd^2(r)\right)\\
&\leq&\sd^2(R)\left[\td^2(R)\left(\int_M\cd(r)\right)^2-\int_M\sd^2(r)\right]
\eeQ
By the H\"older inequality, we have
\beQ
||\varphi||_2^2&\leq&\sd^2(R)\left[\td^2(R)\left(\int_M\cd(r)\right)^2-\frac{1}{\HHp}\left(\int_MH\sd(r)\right)^2\right]
\eeQ
Now using $i)$ of Lemma \ref{lem1}, we get
\beQ
||\varphi||_2^2&\leq&\sd^2(R)\left(\int_M\cd(r)\right)^2\left[\td^2(R)-\frac{1}{\HHp}\right]\leq \widetilde{A_1}C,
\eeQ
where $\widetilde{A_1}$ is a positive constant depending only on the dimension $n$ (because $R\leq\gamma_n$) and $\delta$.$\finpreuve$\\
\indent
The next lemma gives an upper bound for $||\psi||_2$ under the pinching condition.
\begin{lem}\label{lempsi2}
The pinching condition $(\widetilde{P}_C)$ with $C<1$ implies
$$||\psi||^2_2\leq \widetilde{A_2}C,$$
where $\widetilde{A_2}$ is a positive constant depending only on $n$ and $\delta$.
\end{lem}
{\it Proof:} Since $|Z^T|^2=|Z|^2-\left\langle Z,\nu\right\rangle^2$, we have
\beQ
||\psi||_2^2&\leq&\cd^2(R)\left[ \int_M\sd^2(R)-\frac{1}{\HHp}\left( \int_MH\left\langle Z,\nu\right\rangle\right) ^2\right] \\
&\leq& \cd^2(R)\left[ \td^2(R)-\frac{1}{\HHp}\right]\left( \int_M\cd(r)\right)^2 \\
&\leq&\widetilde{A_2}C.
\eeQ
Since $R\leq\gamma_n$, then $\widetilde{A_2}$ depends only on $n$ and $\delta$.
$\finpreuve$

\section{Proof of the Theorems}
\label{proof}

Let $(M^n,g)\in\Mndr$ and $p_0$ the center of the ball of radius $R$ containing $M$. First, we need  the following three lemmas
\begin{lem}\label{prop1}
For any $\varepsilon>0$, there exists $C_{\varepsilon}$ depending on $n$, $\delta$ and $\Hinf$ so that if $(P_{C_{\varepsilon}})$ $\big($or  $(\widetilde{P}_{C_{\varepsilon}})$ for $\delta\geq0\big)$ is true, then
$$\phi(M)\subset \overline{B}(p_0,R)\setminus B(p_0,R-\varepsilon).$$
Moreover, $C_{\varepsilon}\lgra0$ when $\Hinf\lgra+\infty$ or $\varepsilon\lgra 0$.
\end{lem}
We prove this lemma in Section \ref{lemmas}. The second lemma is due to B. Colbois and J.F. Grosjean (see \cite{CG}).
\begin{lem}\label{lem3}
Let $x_0$ be a point of the sphere $S(0,R)$ of $\R^{n+1}$. Assume that $x_0=Ru$ where $u\in\Ss^n$. Now let $(M^n,g)$ be a compact, connected and oriented $n$-dimensional Riemannian manifold without boundary isometrically immersed by $\phi$ into $\R^{n+1}$ so that $$\phi(M)\subset\Big( B(p_0,R+\eta)\setminus B(p_0,R-\eta)\Big) \setminus B(x_0,\rho)$$ with $\rho=4(2n-1)\eta$ and suppose there exists a point $p\in M$ so that $\left\langle Z,u\right\rangle(p)\geq0$.  Then there exists $y_0\in M$ so that the mean curvature $H(y_0)$ at $y_0$ satisfies $|H(y_0)|>\frac{1}{4n\eta}$. 
\end{lem}
\noindent
{\bfseries Remark.}
Note that in \cite{CG}, it is supposed that $\left\langle Z,u\right\rangle >0$, but the condition $\left\langle Z,u\right\rangle\geq0$ is sufficient.
\\\\
\indent
We give a corresponding lemma for the hyperbolic and spherical cases.
\begin{lem}\label{lem4}
Let $x_0$ be a point of the sphere $S(p_0,R)$ of $\HH^{n+1}(\delta)$ ({\it resp.} an open hemisphere of $\Ss^{n+1}(\delta)$). Let $(M^n,g)\in\Mndr$ so that $$\phi(M)\subset\Big( \overline{B}(p_0,R)\setminus B(p_0,R-\eta)\Big) \setminus B(x_0,\rho)$$
with $\rho$ such that 
$$\td\left( (R+\rho)/2\right) -\td(R/2)=4(2n-1)\eta$$
$$\Big({\it resp.}\ \td(R/2)-\td\left( (R-\rho)/2\right) =4(2n-1)\eta\quad\text{if}\ \delta>0\Big)$$
 Then there exist two constants $D$ and $E$ depending on $n$, $\delta$ and $R$ such that if $\eta\leq D$, then there exists $y_0\in M$ so that the mean curvature $H(y_0)$ satisfies
$$| H(y_0)|\geq \frac{E}{8n\eta}.$$
\end{lem}
We prove this Lemma in Section \ref{lemmas}. Now let us prove Theorems \ref{thm1} and \ref{thm3} using the above three Lemmas.
\subsection{Proof of Theorems \ref{thm1} and \ref{thm3}}
\label{sec:proofthm1}
For $\delta=0$, the proof is an immediate consequence of Lemmas \ref{prop1} and \ref{lem3} and is similar to the proof of Theorem 1.2 in \cite{CG}. For $\delta\neq0$, let $\varepsilon>0$. We set $0<\eta\leq\iinf\left\lbrace D,\varepsilon,\frac{\gamma(\varepsilon)}{8(2n-1)}\right\rbrace$, where

$$\gamma(\varepsilon)=
\left\lbrace \begin{array}{cc}
\td\left( \frac{R+\varepsilon}{2}\right)-\td\left( \frac{R}{2}\right)&\text{if}\ \delta<0\\\\
\td\left( \frac{R}{2}\right)-\td\left( \frac{R-\varepsilon}{2}\right)& \text{if}\ \delta>0.
\end{array}
\right. 
$$
Note that $\gamma$ is an increasing smooth function with $\gamma(0)=0$. From Lemma \ref{prop1}, there exists $K_{\varepsilon}=C_{\eta}$ such that $(P_{K_{\varepsilon}})$ implies $$R-r\leq\eta\leq\varepsilon.$$ That's the first point of Theorems \ref{thm1} and \ref{thm3}. Now let's assume that $\varepsilon<\gamma^{-1}\left( \frac{2E}{3\Hinf}\right)$. Suppose there exists $x\in S(p_0,R)$ such that $B(x,\varepsilon)\cap M=\emptyset$. Since $\gamma(\varepsilon)\geq 4(2n-1)\eta$, by Lemma \ref{lem4}, there exists a point $y_0\in M$ so that
$$|H(y_0)|\geq\frac{E}{8n\eta}\geq\frac{(2n-1)E}{n\gamma(\varepsilon)}>\Hinf.$$
Hence a contradiction and $B(x,\varepsilon)\cap M\neq\emptyset$. Moreover, by Lemma \ref{prop1}, $K_{\varepsilon}\lgra0$ when $\Hinf\lgra0$ or $\varepsilon\lgra0$.
$\finpreuve$\\
\indent
From Lemma \ref{prop1}, for any $\varepsilon>0$, there exists $C_{\varepsilon}$ depending on $n$, $\delta$ and $\Hinf$ so that if $(P_{C_{\varepsilon}})$ is true then $|R-r|\leq\varepsilon$. Since $\alpha_n\leq\Hinf\leq\frac{1}{\sqrt{n}}\Binf$, we can assume that $C_{\varepsilon}$ depends on $n$, $\delta$ and $\Binf$.\\
\indent
For the proof of Theorems \ref{thm2} and \ref{thm4} we need the following lemma (which will be proved in Section \ref{lemmas}) on the $L^{\infty}$-norm of $\psi$.
\begin{lem}\label{lem2}
 For any $\varepsilon>0$, there exists $C_{\varepsilon}$ depending on $n$, $\delta$ and $\Binf$ so that if $(P_{C_{\varepsilon}})$  $\big($or  $(\widetilde{P}_{C_{\varepsilon}})$ for $\delta\geq0\big)$is true, then
 $$||\psi||_{\infty}\leq\varepsilon.$$
 Moreover, $C_{\varepsilon}\lgra0$ when $\Binf\lgra+\infty$ or $\varepsilon\lgra 0$.
\end{lem}
\subsection{Proof of Theorems \ref{thm2} and \ref{thm4}}
\label{sec:proofthm2}
Let $\varepsilon>0$ such that $\varepsilon<\td^{-1}\left( \frac{1}{\Hinf}\right) <R$. This choice of $\varepsilon$ implies that if $(P_{C_{\varepsilon}})$ $\big($or  $(\widetilde{P}_{C_{\varepsilon}})\big)$ is true, then $r(x)$ never vanishes. So we can consider the following map
\function{F}{M}{S(p_0,R(M))}{x}{\expo\left(R\left( \diff\expo\right) ^{-1}(\nablab r)\right) .} 
Let $X=\expo^{-1}(x)$. We can easily see that
$$\diff\expo{_{\big|X}}(X)=|X|\nablab r=r\nablab r.$$
\indent
In the case of the Euclidean space $(\delta=0)$, $F(x)$ is precisely $R\frac{X}{|X|}$ where $X$ is the position vector.\\
\indent
We will prove that $F$ is a quasi-isometry. Indeed, we will prove that for any $\theta\in]0,1[$, there exists $\varepsilon(\theta)$ depending on $n$, $\delta$, $\Binf$ and $\theta$ such that for any $x\in M$ and any unit vector $u\in T_xM$, the pinching condition $(P_{C_{\varepsilon(\theta)}})$ implies
$$\Big||\diff_xF(u)|^2-1\Big|\leq\theta.$$
For this, we compute $\diff_xF(u)$ for a unit vector $u\in T_xM$. We have
\beq\label{df1}
\diff F_x(u)&=&\diff\expo_{\Big|R\frac{X}{|X|}}\left( R\;\diff\left(\frac{X}{|X|} \right)_{\big|x}(u) \right) 
\eeq
Let $L(x)=\dfrac{X}{|X|}=\dfrac{\expo^{-1}(x)}{r}$. So we have
\beqt\label{df2}
\diff L_x(u)=\frac{1}{r}\diff\expo^{-1}_{\big|x}(u)-\frac{\diff r(u)}{r^2}\exp^{-1}_{p_0}(x).
\eeqt
Using (\ref{df1}) and (\ref{df2}), we get
\beqt\label{df3}
\diff F_x(u)=\frac{R}{r}\diff\exp_{\Big|R\frac{X}{|X|}}\left( \diff\expo^{-1}_{\big|x}(u)\right)-\frac{R}{r}\diff r(u)\nablab r_{\big|F(x)}
\eeqt
We now compute $\big|\diff_xF(u)\big|^2$. By (\ref{df3}) and the fact that $\expo$ is a radial isometry (see relation (\ref{radialisometry})), we have
\beqt\label{df4}
\big|\diff_xF(u)\big|^2=\frac{R^2}{r^2}\left[ \Bigg|\diff\expo_{\Big|R\frac{X}{|X|}}\left( \diff\exp^{-1}_{p_0\big|x}(u)\right)\Bigg|^2 -\diff r(u)^2\right]
\eeqt
Let $v=u-\left\langle u,\nabla r\right\rangle \nablab r$. That is, $v$ is the part of $u$ normal to $\nablab r$. A straightforward calculation using (\ref{radialisometry}) and (\ref{df4}) shows that
\beqt\label{df5}
\big|\diff_xF(u)\big|^2=\frac{R^2}{r^2}\Bigg|\diff\expo_{\Big|R\frac{X}{|X|}}\left( \diff\expo^{-1}_{\big|x}(v)\right)\Bigg|^2
\eeqt
Finally, by (\ref{expequ1}), we have 
$$\big|\diff\expo^{-1}(v)\big|^2=|v|^2\frac{r^2}{\sd^2(r)},$$
and by (\ref{expequ1}) again,
$$\big|\diff_xF(u)\big|^2=|v|^2\frac{\sd^2(R)}{\sd^2(r)}.$$
From now on, we consider the case $\delta<0$, but the rest of the proof is similar to the case $\delta\geq0$.\\
Since $|v|^2=1-\left\langle u,\nabla r\right\rangle^2 \geq1-|\nabla r|^2$, we deduce that
\beQ
 \Big|\big|\diff_xF(u)\big|^2-1\Big|&\leq&\Big|\frac{\sd^2(R)}{\sd^2(r)}-1\Big|+|\nabla r|\frac{\sd^2(R)}{\sd^2(r)}\\
 &\leq&\sd^2(r)\big|\sd^2(R)-\sd^2(r)\big|+\frac{\sd^2(R)}{\cd(r)\sd^3(r)}||\psi||_{\infty}.
\eeQ
 \indent
 From Lemma \ref{lem3}, we know that for any $\eta>0$, there exists a constant $K_{\eta}$ so that $(P_{K_\eta})$ implies $||\psi||_{\infty}\leq\eta$. Moreover, since $C_{\varepsilon}\lgra0$ when $\varepsilon\lgra0$, there exists $\varepsilon\leq\eta$ depending on $n$, $\delta$, $\Hinf$ and $\eta$ so that $C_{\varepsilon}\leq K_{\eta}$, and then $(P_{C_{\varepsilon}})$ implies $||\psi||_{\infty}\leq\eta$. On the other hand, we have seen that $R$ is bounded by a constant depending only on $n$, $\delta$ and $\Hinf$, then there exist three constants $A_4$, $A_5$ and $A_6$ depending on $n$, $\delta$ and $\Hinf$ so that 
 \beQ
 \Big|\big|\diff_xF(u)\big|^2-1\Big|&\leq&A_4||R-r||_{\infty}+A_5||\psi||_{\infty}\\
 &\leq&A_4\varepsilon+A_5\eta\leq A_6\eta
 \eeQ
 Now, choosing $\eta=\frac{\theta}{A_6}$, we get
 \beqt\label{df6}\Big|\big|\diff_xF(u)\big|^2-1\Big|\leq\theta.\eeqt
For $\theta\in]0,1[$, by (\ref{df6}), $F$ is a local diffeomorphism from $M$ to $S(p_0,R)$. Since for $n\geq2$, $S(p_0,R)$ is simply connected, $F$ is a diffeomorphism. Moreover, the relation (\ref{df6}) says that $F$ is a quasi-isometry.$\finpreuve$

\section{Proof of the technical lemmas}
\label{lemmas}
The proof of Lemmas \ref{prop1} and \ref{lem2} is based on the following Proposition given by a Nirenberg-Moser's type argument.
\begin{prop}\label{fondlem}
Let $(M^n,g)$ be a compact, connected, oriented $n$-dimensional Riemannian manifold without boundary isometrically immersed by $\phi$ into $\R^{n+1}$, $\HH^{n+1}(\delta)$ or an open hemisphere of $\Ss^{n+1}(\delta)$. Let $\xi$ be a nonnegative continuous function on $M$ such that $\xi^k$ is smooth for $k\geq2$. Let $0\leq l<m\leq2$ such that
$$\frac{1}{2}\xi^{2k-2}\Delta\xi^2\leq\ddiv\omega+(\alpha_1+k\alpha_2)\xi^{2k-l}+(\beta_1+k\beta_2)\xi^{2k-m},$$
where $\omega$ is a 1-form and $\alpha_1$, $\alpha_2$, $\beta_1$, $\beta_2$ some nonnegative constants. Then for all $\eta>0$, there exists a constant $L$ depending only on $\alpha_1$, $\alpha_2$, $\beta_1$, $\beta_2$, $\Hinf$ and $\eta$ such that if $||\xi||_{\infty}>\eta$ then
$$||\xi||_{\infty}\leq L||\xi||_2.$$
Moreover, $L$ is bounded when $\eta\longrightarrow\infty$ and if $\beta_1>0$, $L\longrightarrow\infty$ when $\Hinf\longrightarrow\infty$ or $\eta\longrightarrow0$.
\end{prop}

This Proposition is proved in \cite{CG} in the Euclidean case. The proof in the hyperbolic and spherical cases is analogous, using the Sobolev inequality for hypersurfaces of $\HH^{n+1}(\delta)$. Note that in the spherical case, $M$ is assumed to be contained in an open ball of $\Ss^{n+1}(\delta)$ of radius less than $\frac{\pi}{2\sqrt{\delta}}$ as said in Section \ref{sphericalcase} (see \cite{HS1} and \cite{HS2} for details).
\subsection*{Proof of Lemma \ref{prop1}}
We give the proof in the case $\delta\leq0$. For the case $\delta\geq0$ the same computations with $\varphi$ give the result.
First, we compute $\varphi^{2k-2}\Delta\varphi^2$.
\beq\label{Deltaphi}
\varphi^{2k-2}\Delta\varphi^2&=&2\varphi^{2k-1}\Delta\varphi-2|\nabla\varphi|^2\varphi^{2k-2}\nonumber\\
&=&-4\td(r)(1+\delta\td^2(r))\varphi^{2k-1}\Delta r-2|\nabla\varphi|^2\varphi^{2k-2}\nonumber\\
&&+4\varphi^{2k-1}\big|\nabla r\big|^2\Big( 1+3\delta\td^2(r)+2\delta^2\td^4(r) \Big) 
\eeq
Let's compute
\beq\label{div}
(1+\delta\td^2)\td\varphi^{2k-1}\Delta r&=&-\ddiv\left( \varphi^{2k-1}\td(1+\delta\td^2)\nabla r\right) \\
&&+\left\langle \nabla r,\nabla \left( \varphi^{2k-1}\td(1+\delta\td^2)\right) \right\rangle. \nonumber
\eeq
Since $0\leq\td(r)\leq \td(R)<\frac{1}{\sqrt{-\delta}}$ and $|\nabla r|\leq1$, we deduce from the relations(\ref{Deltaphi}) and (\ref{div}) that
\beqt
\varphi^{2k-2}\Delta\varphi^2\leq\ddiv(\omega)+(\alpha_1+k\alpha_2)\varphi^{2k-1}+(\beta_1+k\beta_2)\varphi^{2k-2},
\eeqt
where $\omega$ is a 1-form, $\alpha_1$, $\alpha_2$, $\beta_1$ and $\beta_2$ some nonnegative constants. We can apply Proposition \ref{fondlem} to the function $\varphi$ with $l=1$ and $m=2$. We deduce that if $||\varphi||_{\infty}>\varepsilon$ then there exists a constant $L$ such that $$||\varphi||_{\infty}\leq L||\varphi||_2.$$
On the other hand, by Lemma \ref{lemphi}, we know that if the pinching condition $(P_C)$ is satisfied for $C\leq\alpha(\Hinf)$, then 
$$||\varphi||_2^2\leq A_1C.$$
Take $C=C_{\varepsilon}=\iinf\left\lbrace \alpha(\Hinf),\frac{\varepsilon^2}{L^2A_1}\right\rbrace$. This choice implies
$$ ||\varphi||_{\infty}\leq\varepsilon,$$
that is, $\td^2(R)-\td^2(r)\leq\varepsilon$. Finally, we can choose $C_{\varepsilon}$ smaller in order to have
$R-r\leq\varepsilon.$
$\finpreuve$

\subsection*{Proof of Lemma \ref{lem2}}
We recall that $\psi=\cd(r)|Z^T|$ where $Z^T=\sd(r)\nabla r$. Note that
\beqt\label{nablaZ}
\nabla|Z|^2=2\cd(r)Z^T.
\eeqt
By using th Bochner formula, we the deduce that
\beQ
\frac{1}{2}\Delta\psi^2&=&\frac{1}{2}\Delta\Big|\nabla|Z|^2\Big|^2\\
&=&\left\langle \nabla^*\nabla \diff|Z|^2,\diff|Z|^2\right\rangle-\Big|\nabla \diff|Z|^2\Big|^2\\
&=&\left\langle \Delta \diff|Z|^2,\diff|Z|^2\right\rangle -\Ric\left( \diff|Z|^2,\diff|Z|^2\right) -\Big|\nabla \diff|Z|^2\Big|^2\\
&\leq&\left\langle \Delta \diff|Z|^2,\diff|Z|^2\right\rangle -\Ric\left( \diff|Z|^2,\diff|Z|^2\right)
\eeQ
Now, with the Gauss formula, we can express the Ricci curvature in terms of the second fondamental form $B$. Precisely, we have
\beQ
\frac{1}{2}\Delta\psi^2&\leq&\left\langle \Delta \diff|Z|^2,\diff|Z|^2\right\rangle-nH\left\langle B\nabla|Z|^2,\nabla|Z|^2 \right\rangle \\
&&+\Big| B\nabla|Z|^2\Big|^2-4\delta\Big|\nabla|Z|^2\Big|^2\\\\
&\leq&\left\langle \Delta \diff|Z|^2,\diff|Z|^2\right\rangle-4nH\cd^2(r)\left\langle BZ^T,Z^T\right\rangle\\
&& +4\cd^2(r)|BZ^T|^2-4(n-1)\cd^2(r)\delta|Z^T|^2
\eeQ
Since $|Z^T|\leq\sd(R)$, we easily see that the pinching condition $(P_C)$, with $C<1$ for $\delta\geq0$ or $C<\alpha(\Hinf)$ for $\delta<0$, implies
$$||\psi||_{\infty}\leq K_1,$$
where $K_1$ is a positive constant depending only on $n$, $\delta$ and $\Binf$. It follows that
\beQ
\frac{1}{2}\left( \Delta\psi^2\right) \psi^{2k-2}\leq\left\langle \Delta \diff|Z|^2,\diff|Z|^2\right\rangle\psi^{2k-2}+K_2\psi^{2k-2}.
\eeQ
Let $\omega=\Delta|Z|^2\psi^{2k-2}\diff|Z|^2$. We have
\beQ
\left\langle \Delta \diff|Z|^2,\diff|Z|^2\right\rangle\psi^{2k-2}
&=&\ddiv(\omega)+(\Delta|Z|^2)^2\psi^{2k-1}\\
&&-2(2k-2)\cd(r)\Delta|Z|^2\left\langle Z^T,\diff\psi\right\rangle \psi^{2k-3},
\eeQ
Moreover, a straightforward calculation using the facts that  
$$e_i(|Z^T|)=\frac{1}{2}\frac{e_i(|Z^T|^2)}{|Z^T|},$$
and $|Z^T|^2=\sd^2(r)-\left\langle Z,\nu\right\rangle^2$ gives

\beQ
e_i(\psi)&=&\cd(r)\frac{\delta\cd(r)\sd(r)e_i(r) -\left\langle Z,\nu\right\rangle\left(  B_{ij}\left\langle Z,e_j\right\rangle +\left\langle \nablab_{e_i}Z,\nu\right\rangle\right)  }{|Z^T|}\\
&&-\delta\sd(r)|Z^T|e_i(r)
\eeQ
All the terms can be bounded easily except $\left\langle \nablab_{e_i}Z,\nu\right\rangle$ which will be investigated. Since $Z=\sd(r)\nablab r$, this is equivalent to have an upper bound for $\left\langle \nablab_{e_i}\nablab r,\nu\right\rangle$.
From (\ref{expequ2}) and (\ref{nullspace}), we deduce that
$$\big|\left\langle \nablab_{e_i}\nablab r,\nu\right\rangle\big|=\frac{\cd(r)}{\sd(r)}\left\langle e_i^t,\nu^t\right\rangle, $$
where $e_i^t$ and $\nu^t$ are the part of $e_i$ and $\nu$ tangent to the geodesic sphere of radius $r$, that is, orthogonal to $\nablab r$. Since 
$$|\nu^t|^2=1-\left\langle \nu,\nablab r \right\rangle^2=|\nabla r|^2, $$
we have 
$$\big|\left\langle \nablab_{e_i}\nablab r,\nu\right\rangle\big|\leq\frac{\cd(r)}{\sd(r)}|\nabla r|^.$$
Then
\beQ
\left\langle \Delta \diff|Z|^2,\diff|Z|^2\right\rangle\psi^{2k-2}&\leq&\ddiv(\omega)+(\Delta|Z|^2)^2\psi^{2k-1}
\\&&+2(2k-2)\Big|\Delta|Z|^2\Big|\psi^{2k-2}\\
&&+K_3\Big|\Delta|Z|^2\Big|||B||_{\infty}\psi^{2k-2},
\eeQ
where, $K_3$ depends on $n$, $\delta$ and $\Hinf$. Moreover, $\Delta|Z|^2=-2\ddiv\left( \cd(r)Z^T\right)$ and by Lemma \ref{lem1} $i)$, we deduce that there exists a constant $K_4$ depending on $n$, $\delta$ and $\Hinf$ such that
$$\Delta|Z|^2\leq K_4.$$
Finally, we have
\beQ
\psi^{2k-2}\Delta\psi^2&\leq&\ddiv(\omega)+(\alpha_3+k\alpha_4)\psi^{2k-1}+(\beta_3+k\beta_4)\psi^{2k-2}
\eeQ
with some nonnegative constants $\alpha_3$,$\alpha_4$,$\beta_3$ and $\beta_4$ depending on $n$, $\delta$ and $\Binf$.
Now applying Proposition \ref{fondlem} with $l=1$ and $m=2$, we get that for $\eta>0$, there exists $L$ depending on $n$, $\delta$, $\Binf$ and $\eta$ so that if $||\psi||_{\infty}>\eta$ then
$$||\psi||_{\infty}\leq L||\psi||_2.$$
From Lemma \ref{lempsi} we deduce that if $(P_C)$ holds with $C<1$ for $\delta\geq0$ or $C<\alpha(\Hinf)$ for $\delta<0$, then 
$$||\psi||_2^2\leq A_2C+A_3||\varphi||_{\infty}.$$
 Let $\varepsilon>0$, and put $K_{\varepsilon}:=\iinf\left\lbrace \frac{\varepsilon^2}{2L^2A_2},C_{\frac{\varepsilon^2}{2L^2A_3}}\right\rbrace $ where $C$ is the constant defined in the proof of the Lemma \ref{prop1}. Then, if $(P_{K_{\varepsilon}})$ holds, we have
 $$||\psi||_{\infty}\leq L||\psi||_2\leq\varepsilon.$$
$\finpreuve$

\subsection*{Proof of Lemma \ref{lem4}}
We give the proof in the hyperbolic case, $\delta=-1$. The proof for any $\delta<0$ or for the spherical case is similar.\\
\indent
Let us consider $f$ the conformal map from the unit ball $\widetilde{B}(0,1)$ of $\R^{n+1}$ into $\HH^{n+1}$ so that $f(0)=p_0$. The conformal factor is the function $h^2$ where $h$ is defined by
\function{h}{[0,1[}{\R_+}{r}{\frac{2}{1-r^2}}
\indent
For any $\rho>0$, $B(p_0,\rho)=f\left( \widetilde{B}(0,\widetilde{\rho})\right) $, where $\widetilde{B}(0,\widetilde{\rho})\subset\widetilde{B}(0,1)$ is the ball of radius $a(\rho):=\widetilde{\rho}=\td(\rho/2)$. Let $\widetilde{\phi}=f^{-1}\circ\phi$. By hypothesis, 
$$\phi(M)\subset\Big(B(p_0,R+\eta)\setminus B(p_0,R-\eta)\Big)\setminus B(x_0,\rho),$$
with $x_0\in S(p_0,R)$ and $\rho$ chosen so that
$$a(R+\rho)-a(R)=4(2n-1)\eta,$$
then
$$\widetilde{\phi}(M)\subset\Big(B\big(p_0,a(R+\eta)\big)\setminus B\big(p_0,a(R-\eta0\big)\Big)\setminus B(z_0,\rho'),$$

where $z_0=\frac{1}{2}\left[ a(R+\rho)+a(R-\rho)\right]u$, with $u$ a unit vector and 
$$\rho'=\frac{1}{2}\left[ a(R+\rho)-a(R-\rho)\right].$$
Obviously we have
$$  \widetilde{B}\big( 0,a(R+\eta)\big) \setminus\widetilde{B}\big( 0,a(R-\eta)\big)\subset \widetilde{B}\big( 0,a(R)+\eta\big) \setminus\widetilde{B}\big( 0,a(R)-\eta\big).$$
Moreover, by concavity of the function $a$, we have
$$a(R+\rho)+a(R-\rho)\leq 2a(R),$$
and then
$$\widetilde{B}(z_0,\rho')\supset\widetilde{B}(\widetilde{x_0},\rho''),$$
where $\widetilde{x_0}=a(R)u$ and $\rho''=a(R+\rho)-a(R)$. Finally, we have
$$\widetilde{\phi}(M)\subset\left[ \widetilde{B}\big( 0,a(R)+\eta\big) \setminus\widetilde{B}\big( 0,a(R)-\eta\big)\right] \setminus\widetilde{B}\big(\widetilde{x_0},4(2n-1)\eta\big).$$
\begin{center}
\includegraphics[width=7cm]{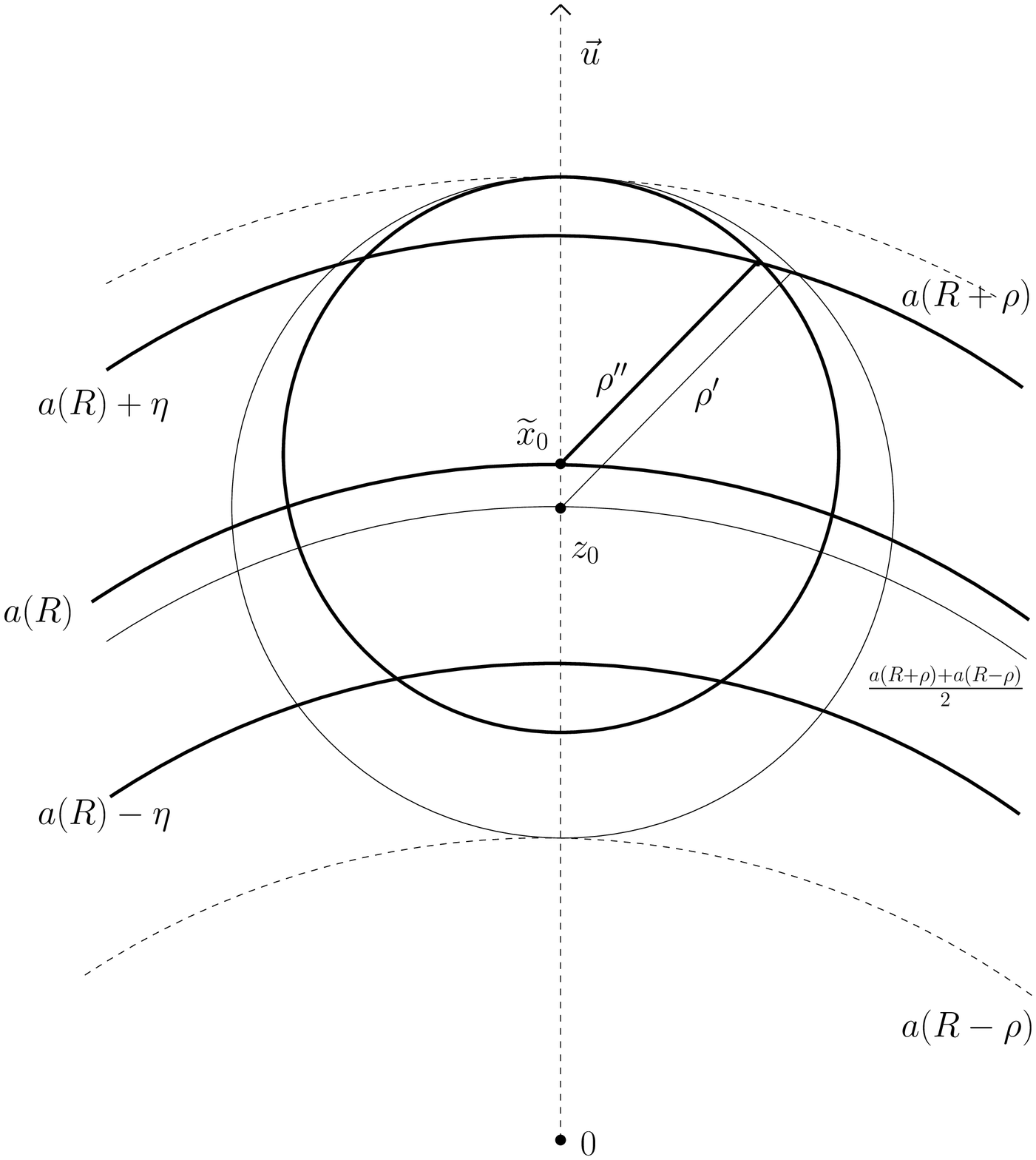}
\end{center}

\noindent
Since $a(R)$ is the extrinsic radius of $\widetilde{\phi}(M)$, there exists a point $p\in\widetilde{\phi}(M)$ so that $\left\langle \widetilde{Z},u\right\rangle(p)\geq0$, where $\widetilde{Z}$ is the position vector of $\widetilde{\phi}(M)$ in $\widetilde{B}(0,1)$. By Lemma \ref{lem3}, there exists $y_0\in\widetilde{\phi}(M)$ so that the mean curvature $\widetilde{H}$ satisfies $|\widetilde{H}(y_0)|>\frac{1}{4n\eta}$. Moreover, we have the well-known formula for the conformal mean curvature (see for example \cite{Esc2})
$$H=h^{-1}\left( \widetilde{H}+h^{-1}\left\langle \widetilde{\nabla} h,\widetilde{\nu}\right\rangle \right) ,$$
where $\widetilde{\nabla}$ and $\widetilde{\nu}$ are the gradient and the normal unit vector field in $\widetilde{B}(0,1)$, and $\left\langle .,.\right\rangle $ is the Euclidean scalar product in $\widetilde{B}(0,1)$. Therefore,
$$|H|\geq h^{-1}(\tilde{r})\left(\widetilde{H}-h^{-1}(\tilde{r})|\tilde{\nabla} h| \right),$$
where $\widetilde{r}(x)$ is the Euclidean distance from $0$ to $x$. So we have
\beqt \label{H1}
\frac{1}{2}\geq h^{-1}( \tilde{r})=\frac{1-\tilde{r}^2}{2}\geq\frac{1-a(R)^2}{2},
\eeqt
and 
\beqt\label{H2}
|\widetilde{\nabla} h|=\frac{4}{\left( 1-\tilde{r}^2\right)^2 }|\tilde{r}\widetilde{\nabla}\tilde{r}|\leq\frac{4}{\left( 1-a(R)^2\right)^2}.
\eeqt
Finally, by (\ref{H1}) and (\ref{H2}), we get
$$|H(f^{-1}(y_0))|\geq \frac{E}{4n\eta}-\frac{1}{4E^2},$$
where $E=\frac{1-a(R)^2}{2}$ is a constant depending on $n$, $\delta$ and $R$. Moreover, there exists a constant $D$ depending on $n$, $\delta$ and $R$ so that if $\eta\leq D$, then $\frac{1}{4E^2}\leq\frac{E}{8n\eta}$, and so
$$|H(f(y_0))|\geq \frac{E}{8n\eta}.$$
$\finpreuve$
\noindent
{\bfseries Remark.}
If we suppose $(P_C)$ with $C<\alpha(\Hinf)$ for $\delta<0$ ({\it resp.} with $C<1$ for $\delta\geq0$) then $D$ and $E$ depend on $n$, $\delta$ and $\Hinf$ ({\it resp.} on $n$ and $\delta$).
\\
{\bfseries Remark.}
For $\delta>0$, the function $a$ is convex and so $\rho''=a(R)-a(R-\rho)$.


\end{document}